\documentstyle[12pt,twoside]{article}
\font\teneufm=eufm10
\font\seveneufm=eufm7
\font\fiveeufm=eufm5
\newfam\eufmfam
\textfont\eufmfam=\teneufm
\scriptfont\eufmfam=\seveneufm
\scriptscriptfont\eufmfam=\fiveeufm

\newfam\msbfam
\font\tenmsb=msbm10 scaled \magstep1  \textfont\msbfam=\tenmsb
\font\sevenmsb=msbm7 scaled \magstep1 \scriptfont\msbfam=\sevenmsb
\font\fivemsb=msbm5 scaled \magstep1  \scriptscriptfont\msbfam=\fivemsb
\def\Bbb{\fam\msbfam \tenmsb}

\def\RR{{\Bbb R}}
\def\CC{{\Bbb C}}

\def\NN{{\Bbb N}}
\def\ZZ{{\Bbb Z}}

\def\Aut{\hbox{Aut}}

\def\ra{\rightarrow}

\def\HollowBoxx #1#2#3{{\dimen0=#1 \advance\dimen0 by -#2       
       \dimen1=#1 \advance\dimen1 by #3                       
        \vrule height 0pt depth #3 width #2                   
       \hskip -#3
       \vrule height #1 depth #3 width #3}}                   
 \def\LeftContraction{\mathord{\kern1.45pt \HollowBoxx{6pt}{3.5pt}{.4pt}}\,}

 \def\HollowBox #1#2#3{{\dimen0=#1 \advance\dimen0 by -#3       
       \dimen1=#1 \advance\dimen1 by #3                       
        \vrule height #1 depth #3 width #3                    
        \vrule height 0pt depth #3 width #2                   
        \hskip -#3}}                                             
 \def\RightContraction{\mathord{\, \HollowBox{6pt}{3.1pt}{.4pt}} \kern1.6pt}             

\def\qed{{\hfill $\Box$}}
\def\Aut{\hbox{Aut}}

\newtheorem{theorem}{THEOREM}[section]

\newtheorem{corollary}[theorem]{COROLLARY}
\newtheorem{lemma}[theorem]{Lemma}

\newtheorem{proposition}[theorem]{PROPOSITION}

\begin{document}
\begin{center}
{\Large \bf On the Dimensions
\medskip\\
of the Automorphism Groups} 
\medskip \\
{\Large \bf of Hyperbolic Reinhardt Domains}  
\medskip \\
{\normalsize James A. Gifford, \ \ \  Alexander V. Isaev \ \ \ and \ \ \
Steven G. Krantz\footnote{Work
supported in part by NSF Grant DMS-9531967}}
\end{center} 

\begin{quotation} 
\small \sl We study the possible dimensions that the groups of
holomorphic automorphisms of hyperbolic Reinhardt domains can have. We are
particularly interested in the problem of characterizing Reinhardt domains with
automorphism group of prescribed dimension.
\end{quotation}

\pagestyle{myheadings}
\markboth{J. A. Gifford, A. V. Isaev and S. G. Krantz}{ON THE
DIMENSIONS OF AUTOMORPHISM GROUPS}

\setcounter{section}{-1}

\section{Introduction}
\setcounter{equation}{0}

Let $D$ be a domain (a connected open set) in $\CC^n$, $n\ge 2$. Denote by $\hbox{Aut}(D)$
the group of holomorphic automorphisms of $D$; that is, $\Aut(D)$ is
the group under composition of all biholomorphic self-maps of $D$. If $D$ is bounded or,
more generally, Kobayashi-hyperbolic, then the group $\hbox{Aut}(D)$
with the topology of uniform convergence on compact subsets of $D$ is
in fact a finite-dimensional Lie group (see \cite{Ko}). We note that this Lie group is always a {\it real} group but
never a complex Lie group (except for the case of zero-dimensional groups).  Thus, when we specify the dimension
of this group, we shall always be speaking of its real dimension.
By contrast, when we speak of the dimension of the domain on which
it acts, we shall be referring to complex dimension.

We are interested in characterizing a domain by its automorphism
group. Much work has been done on classifying domains with non-compact
automorphism group (see \cite{IK2} for a detailed exposition). In this
paper we concentrate on the dimension of $\hbox{Aut}(D)$. Namely, we
are interested in the following question: to what extent does the dimension
of the automorphism group determine the domain?

In this work we consider only Reinhardt domains, i.e., domains
invariant under the (coordinate) rotations
$$
z_j\mapsto e^{i\phi_j}z_j,\qquad \phi_j\in\RR,\qquad j=1,\dots,n.
$$
As we shall see below, even this special case leads to difficult
problems.

The paper is based on the structure theorem by Kruzhilin \cite{Kr}
(see also \cite{Sh} for the case of bounded domains) that allows us to
list all possible dimensions that the automorphism groups of
hyperbolic Reinhardt domains can have. It turns out that all
dimensions (except for the case of domains that, up to dilations and
permutations of coordinates, are the unit ball in $\CC^n$---see
Corollary \ref{c2}) lie
between $n$ and $n^2+2$ inclusive; the dimensions are
even if $n$ is even and odd if $n$ is odd. 

We classify all domains
whose automorphism groups have dimensions $n^2$ and $n^2+2$ (see
Theorem \ref{t9} and Corollary \ref{c5}).
The remaining dimensions (i.e., those that lie between $n$ and $n^2-2$ inclusive)
split into two sets: the ``bad'' and
``good'' ones (the latter corresponds to the case of domains with
non-compact automorphism group). These will be defined in the sequel. 

While it will turn out that there is no hope to obtain
any reasonable classification of domains whose automorphism groups
have dimensions that belong to the ``bad'' set, one {\it can} hope to
obtain some description for the ``good'' dimensions. For
$C^1$-smoothly bounded domains such a description has been already found
in \cite{IK1}.

In this paper we study the structure of the sets of ``bad'' and
``good'' dimensions. The main question we are interested in is:
what is the asymptotic behavior (as a function of $n$) of the numbers
of ``bad'' and ``good'' dimensions as the spatial dimension
$n\rightarrow\infty$? We have been able to prove that the number of ``bad'' dimensions behaves
asymptotically as $n^2/2$ (see Theorem \ref{num}) which means that these dimensions
asymptotically fill the whole list of all possible automorphism group
dimensions. On the other hand, we show that the $\liminf$ of the number of ``good'' dimensions
behaves asymptotically at least as $n$ (see Theorem \ref{numh}). This
last result implies that the
probability of randomly choosing a ``good'' dimension from
the list of all possible automorphism group dimensions is
asymptotically at least of order
$1/n$; this information is encouraging compared with the fact that
almost any randomly chosen domain in $\CC^n$ does not belong to any
reasonable classification list. 

We have also
made numerical computations for the numbers of ``bad'' and ``good''
dimensions for up to $n=4065$ and present some of the results in
Section 2. These results were obtained by combining C programming with
techniques of network programming. The source code of the C program is available on
the World Wide Web at 
\smallskip \\

\begin{center}
{\large \tt http://wwwmaths.anu.edu.au/$\sim$james/reinhardt} 
\end{center}
\vspace*{.12in}

As we shall see below, finding the numbers of ``bad'' and ``good''
dimensions is also related to determining certain
characteristics of partitions by way of their Young diagrams.

We wish to thank G. Andrews, A. Molev, M. F. Newman and R. Stanley for useful
discussions and interest in our work.

\section{Results}
\setcounter{equation}{0}

Let $\CC^*:=\CC\setminus\{0\}$. We denote by $\hbox{Aut}_{alg}((
\CC^*)^n)$ the group of algebraic automorphisms of
$(\CC^*)^n$, i.e., the group of mappings of the form
\begin{equation}\label{a}
z_i\mapsto\lambda_iz_1^{a_{i1}}\dots
z_n^{a_{in}},\quad i=1,\dots n,
\end{equation}
where $\lambda_i\in{\CC}^{*}$, $a_{ij}\in{\ZZ}$, and
$\hbox{det}(a_{ij})=\pm 1$.

For a hyperbolic Reinhardt domain $D\subset{\CC}^n$,
denote by $\hbox{Aut}_{alg}(D)$ the subgroup of
$\hbox{Aut}(D)$ that consists of algebraic automorphisms of $D$,
i.e., automorphisms induced by elements of
$\hbox{Aut}_{alg}(({\CC}^*)^n)$. Let  $\hbox{Aut}_0(D)$ be
the connected component of the identity in $\hbox{Aut}(D)$, and
the dot (the symbol ${}\cdot{}$) denote the composition operation in
$\hbox{Aut}(D)$. It is shown in
\cite{Kr} that
$\hbox{Aut}(D)=\hbox{Aut}_0(D)\cdot \hbox{Aut}_{alg}(D)$.

By \cite{Kr} any hyperbolic
Reinhardt domain in ${\CC}^n$ can---by a biholomorphic mapping of the
form (\ref{a})---be put into a normalized form $G$ written as follows.
There exist integers $0\le s\le
t\le p\le n$ and $n_i\ge 1$, $i=1,\dots,p$, with $\sum_{i=1}^p
n_i=n$, and real numbers ${\alpha}_i^j$, $i=1,\dots,s$,
$j=t+1,\dots,p$, ${\beta}_j^k$, $j=s+1,\dots,t$,
$k=t+1,\dots,p$, such that if we set
$z^i:=\left(z_{n_1+\dots+n_{i-1}+1},\dots,z_{n_1+\dots+n_i}\right)$,
$i=1,\dots,p$, then $G$ can be
written in the form
\begin{eqnarray}
G&=&\Biggl\{\left|z^1\right|<1,\dots,\left|z^s\right|<1,\nonumber\\
&{}&\Biggl(\frac{z^{t+1}}{\prod_{i=1}^s
\left(1-\left|z^i\right|^2\right)^{{\alpha}_i^{t+1}}\prod_{j=s+1}^t
\exp\left(-{\beta}_j^{t+1}\left|z^j\right|^2\right)}\ ,\ \dots \ , \label{b}\\
&{}&\frac{z^{p}}{\prod_{i=1}^s
\left(1-\left|z^i\right|^2\right)^{{\alpha}_i^p}\prod_{j=s+1}^t
\exp\left(-{\beta}_j^p\left|z^j\right|^2\right)}\Biggr)\in\widetilde
G\Biggr\},\nonumber
\end{eqnarray}
where $\widetilde G:=G\bigcap\left\{z^i=0,\,
i=1,\dots,t\right\}$ is some hyperbolic Reinhardt domain in
${\CC}^{n_{t+1}}\times\dots\times{\CC}^{n_p}$. It should be noted that
any given domain will have many different normalized forms of type (\ref{b}).

A normalized form can be chosen so that
$\hbox{Aut}_0(G)$ is given by the following formulas:
\begin{eqnarray}
z^i&\mapsto&\frac{A^iz^i+b^i}{c^iz^i+d^i},\quad i=1,\dots,s,\nonumber\\
z^j&\mapsto& B^jz^j+e^j,\quad j=s+1,\dots,t,\label{c}\\
z^k&\mapsto&
C^k\frac{\prod_{j=s+1}^t\exp\left(-\beta_j^k\left(2\overline{e^j}^TB^jz^j+
|e^j|^2\right)\right)z^k}{\prod_{i=1}^s(c^iz^i+d^i)^{2\alpha_i^k}},\quad
k=t+1,\dots,p,\nonumber
\end{eqnarray}
where
\begin{eqnarray}
&{}&\left(\begin{array}{cc}
A^i& b^i\\
c^i& d^i
\end{array}
\right)
\in SU(n_i,1),\quad i=1,\dots,s,\nonumber\\
&{}&B^j\in U(n_j),\quad e^j\in{\Bbb C}^{n_j},\quad j=s+1,\dots,t,\label{d}\\
&{}&C^k\in U(n_k),\quad k=t+1,\dots,p.\nonumber
\end{eqnarray}

It follows from (\ref{c}), (\ref{d}) that the dimension of
the automorphism group of any hyperbolic Reinhardt domain in $\CC^n$
is a number of the form
\begin{equation}
\sum_{i=1}^k n_i^2+2\sum_{j=1}^m n_j, \qquad 0\le m\le k, \label{e}
\end{equation}
for some partition $(n_1,\dots,n_k)$ of $n$, $1\le k\le n$. We will be
interested in the structure of the set $\Omega(n)$ of all numbers
(\ref{e}). Let $\Omega(n,\ell)$ be the set of all numbers of the form
(\ref{e}) with $k=\ell$ and $\Omega(n,\ell,q)$ the set of all numbers of the
form (\ref{e}) with $k=\ell$, $m=q$. Clearly,
$\Omega(n,\ell)=\cup_{q=0}^\ell\Omega(n,\ell,q)$, $\Omega(n)=\cup_{\ell=1}^{n}\Omega(n,\ell)$. First,
prove the following:

\begin{proposition}\label{p1} \sl Let $N\in\Omega(n)$. Then 
\smallskip\\
\noindent {\bf (i)} $N$ is even (odd) if $n$ is even (odd).
\smallskip\\

\noindent {\bf (ii)} $N\ge n$.
\smallskip\\

\noindent {\bf (iii)} If $N\in\Omega(n,\ell)$ for $\ell\ge 2$  then $N\le n^2+2$.
\end{proposition}

\noindent {\bf Proof:} Statements {\bf (i)} and {\bf (ii)} are
obvious. We prove {\bf (iii)} by induction. It is obvious
for $n=2$, so we assume that $n\ge 3$. Let $N=\sum_{i=1}^\ell
n_i^2+2\sum_{j=1}^m n_j$, for some $2\le \ell\le n$, $0\le m\le \ell$.

Suppose first that $\ell=2$, $m\le1$. Then we have
\begin{eqnarray*}
&{}&N=\left(n_1^2+2\sum_{j=1}^m n_j\right)+n_2^2\le\\
&{}&(n-n_2)^2+2(n-n_2)+n_2^2=n^2+2(n_2^2-n_2(n+1)+n)<n^2+2.
\end{eqnarray*}

Let $\ell=2$, $m=2$. Then
$$
N=(n-n_2)^2+n_2^2+2n=n^2+2(n_2^2-n_2n+n)\le n^2+2.
$$

Suppose now that $\ell\ge3$ and assume first that $m\le \ell-1$. Then by induction
we have
\begin{eqnarray*}
&{}&N=\left(\sum_{i=1}^{\ell-1}n_i^2+2\sum_{j=1}^m n_j\right)+n_{\ell}^2\le\\
&{}&(n-n_{\ell})^2+2+n_{\ell}^2=n^2+2+2n_{\ell}(n_{\ell}-n)<n^2+2.
\end{eqnarray*}

Assume finally that $\ell\ge3$ and $m=\ell$. Then by induction
we get
\begin{eqnarray*}
&{}&N=\left(\sum_{i=1}^{\ell-1}n_i^2+2\sum_{j=1}^{\ell-1} n_j\right)+n_{\ell}^2+2n_{\ell}\le\\
&{}&(n-n_{\ell})^2+2+n_{\ell}^2+2n_{\ell}=n^2+2+2n_{\ell}(n_{\ell}-n+1)\le n^2+2.
\end{eqnarray*} 

The proposition is proved. \qed
\smallskip\\

It follows from Proposition \ref{p1} that the value $n^2+2n$ can only be
taken by
$\Omega(n,1,1)$ (which is clearly a one-point set) corresponding to the case $s=p=1$ in formula
(\ref{b}). Thus we obtain the
following characterization of the unit ball in the class of hyperbolic
Reinhardt domains.

\begin{corollary}\label{c2} \sl Let $D\subset\CC^n$ is a hyperbolic Reinhardt domain
such that $\hbox{dim}\,\hbox{Aut}(D)>n^2+2$. Then, up to dilations and
permutations of coordinates, $D$ is the unit ball $B^n$.
\end{corollary}

In the following proposition we establish upper bounds for\linebreak 
$\Omega(n,\ell,0)$, $\ell\ge 2$.

\begin{proposition}\label{p3} \sl Let $N\in\Omega(n,\ell,0)$. Then $N\le
(n-\ell+1)^2+\ell-1$.
\end{proposition}

\noindent{\bf Proof:} The inequality is obvious for $\ell=1$, so we
assume that $\ell\ge 2$. The proof for $\ell\ge 2$ proceeds by induction. The
inequality is clearly correct for $n=2$ and we suppose that
$n\ge3$. Let $N=\sum_{i=1}^{\ell} n_i^2$. If $\ell=2$, we have
$$
N=n_1^2+n_2^2=(n-n_2)^2+n_2^2=(n-1)^2+1+2(n_2^2-n_2n+n-1)\le(n-1)^2+1.
$$

Assume now that $\ell\ge3$. Then by induction we have
\begin{eqnarray*}
&{}&N=\sum_{i=1}^{\ell-1}n_i^2+n_{\ell}^2\le\\
&{}&(n-n_{\ell}-\ell+2)^2+\ell-2+n_{\ell}^2=(n-\ell+1)^2+\ell-1+\\
&{}&2(n_{\ell}^2-n_{\ell}(n-\ell+2)+n-\ell+1)\le(n-\ell+1)^2+\ell-1.
\end{eqnarray*}

The proposition is proved. \qed
\smallskip\\

\begin{corollary}\label{c4} \sl Let $N\in\Omega(n,\ell)$. Then $N\le(n-\ell+1)^2+\ell-1+2n$.
\end{corollary}

It follows from Corollary \ref{c4} that the value $n^2+2$ can only be taken
by the elements of $\Omega(n,2,2)$. The only partition $(n_1,n_2)$
(assuming $n_1\le n_2$)
that can realize this value is clearly $(1,n-1)$ which corresponds to
the case $s=p=2$, $n_1=1$, $n_2=n-1$ in formula (\ref{b}). Thus, we
obtain the following:

\begin{corollary}\label{c5}\sl Let $D\subset\CC^n$ be a hyperbolic Reinhardt domain
such that $\hbox{dim}\,\hbox{Aut}(D)=n^2+2$. Then, up to dilations and
permutations of coordinates, $D$ is the product $B^{n-1}\times\Delta$
of the unit ball
$B^{n-1}\subset\CC^{n-1}$ and the unit disc $\Delta\subset\CC$.
\end{corollary}

Now we shall deal with dimension $n^2$ and classify all
hyperbolic Reinhardt domains $D$ such that
$\hbox{dim}\,\hbox{Aut}(D)=n^2$. For this result, we need to understand what
numbers from $\Omega(n)$ can equal $n^2$. First of all, we
clearly have $\Omega(n,1,0)=\{n^2\}$. Next, the following holds. 

\begin{proposition}\label{p6} \sl If $n\ge 4$, then for any $N\in\Omega(n,\ell)$ with $\ell\ge
3$, one has $N<n^2$.
\end{proposition}

\noindent {\bf Proof:} It follows from Corollary \ref{c4} that
$$
N\le n^2+\ell^2-\ell(2n+1)+4n\le n^2-2(n-3)<n^2.
$$

The proposition is proved.\qed
\smallskip\\

It follows from Proposition \ref{p6} that, for $n\ge 4$, if $n^2\in\Omega(n,\ell)$
and $\ell\ge 2$, then $\ell=2$. 

We now take a closer look at the set $\Omega(n,2)$. To cover all the
elements from  $\Omega(n,2,0)$ we clearly only need to consider
partitions $(n_1,n_2)$ of $n$ with $n_1\ge n_2$, i.e. partitions of
the form
\begin{eqnarray*}
&{}&\left(\frac{n}{2}+\mu,\frac{n}{2}-\mu\right),\,\mu=0,\dots,\frac{n}{2}-1,\,\,\hbox{if
$n$ is even},\\
&{}&\left(\frac{n+1}{2}+\mu,\frac{n-1}{2}-\mu\right),\,\mu=0,\dots,\frac{n-3}{2},\,\,\hbox{if
$n$ is odd}.
\end{eqnarray*}
Therefore, the following proposition is obvious.

\begin{proposition}\label{p7} \sl \hbox{ \ \ } \hfill \break
\vspace*{.003in} 

\noindent {\bf (i)} If $n$ is even,
\begin{eqnarray*}
\Omega(n,2)&=&\left\{\frac{n^2}{2}+2\mu^2:\,
\mu=0,\dots,\frac{n}{2}-1\right\}\cup\\
&{}&\left\{\frac{n^2}{2}+2\mu(\mu-1)+n:\,
\mu=0,\dots,\frac{n}{2}-1\right\}\cup\\
&{}&\left\{\frac{n^2}{2}+2\mu(\mu+1)+n:\,
\mu=0,\dots,\frac{n}{2}-1\right\}\cup\\
&{}&\left\{\frac{n^2}{2}+2\mu^2+2n:\,
\mu=0,\dots,\frac{n}{2}-1\right\}.
\end{eqnarray*}
\smallskip\\

\noindent{\bf (ii)} If $n$ is odd,
\begin{eqnarray*}
\Omega(n,2)&=&\left\{\frac{n^2+1}{2}+2\mu(\mu+1):\,
\mu=0,\dots,\frac{n-3}{2}\right\}\cup\\
&{}&\left\{\frac{n^2+1}{2}+2\mu^2+n-1:\,
\mu=0,\dots,\frac{n-3}{2}\right\}\cup\\
&{}&\left\{\frac{n^2+1}{2}+2\mu(\mu+2)+n+1:\,
\mu=0,\dots,\frac{n-3}{2}\right\}\cup\\
&{}&\left\{\frac{n^2+1}{2}+2\mu(\mu+1)+2n:\,
\mu=0,\dots,\frac{n-3}{2}\right\}.
\end{eqnarray*}
\end{proposition}

\begin{corollary}\label{c8} \sl \hbox{ \ \ } \hfill \break
\vspace*{.003in}

\noindent {\bf (i)} A number $N\in\Omega(n,2)$, $n\ne 4$, is equal to $n^2$ only if
$N\in\Omega(n,2,1)$ and corresponds to the
partition $(n-1,1)$.
\smallskip\\

\noindent {\bf (ii)} A number $N\in\Omega(4,2)$ is equal to 16 only if
either
\begin{itemize}
\item[{\bf (a)}]  $N\in\Omega(4,2,2)$ and $N$ corresponds to the 
partition $(2,2)$; or
\item[{\bf (b)}]  $N\in\Omega(4,2,1)$ and $N$ corresponds to the 
partition $(3,1)$.
\end{itemize} 
\end{corollary}

We are now ready to prove the following classification result for
dimension $n^2$.

\begin{theorem}\label{t9} \sl Let $D\subset\CC^n$ be a hyperbolic Reinhardt domain
such that $\hbox{dim}\,\hbox{Aut}(D)=n^2$. Then $D$ is holomorphically
equivalent to one of the following domains:
\smallskip

\noindent {\bf (i)} $\{z\in\CC^n: r<|z|<R\}$, $0\le r<R<\infty$;
\smallskip\\

\noindent {\bf (ii)} $\Delta^3$ (here $n=3$);
\smallskip\\

\noindent {\bf (iii)} $B^2\times B^2$ (here $n=4$);
\smallskip\\

\noindent {\bf (iv)} $\{(z',z_n)\in\CC^n:
|z'|^2+|z_n|^{\alpha}<1\}$, $\alpha\in\RR$, $\alpha\ne 0,2$;
\smallskip\\

\noindent {\bf (v)}
$\{(z',z_n)\in\CC^n:|z'|<1,r(1-|z'|^2)^{\alpha}<|z_n|<R(1-|z'|^2)^{\alpha}\}$,
$\alpha\in\RR$, $0<r<R\le\infty$;
\smallskip\\

\noindent {\bf (vi)}
$\{(z',z_n)\in\CC^n:re^{\alpha|z'|^2}<|z_n|<Re^{\alpha|z'|^2}\}$, where
$0<r<R\le\infty$, $\alpha\in\RR$, $\alpha\ne 0$ and, if $R=\infty$, then
$\alpha>0$.
\smallskip\\

The equivalence is given by a mapping of the form (\ref{a}).
\end{theorem}

\noindent{\bf Proof:} Let $G$ be a normalized form of $D$ as in (\ref{b}).
We first consider the case of
$\Omega(n,1,0)$. Then $\hbox{Aut}_0(G)=U(n)$ (see (\ref{c})). It is clear that any
hyperbolic domain with this property has the form {\bf (i)}.

Next, by Proposition \ref{p6}, the case of $\Omega(n,3)$ is only non-trivial
when
$n=3$. Clearly, $G$ then coincides with $\Delta^3$.

We now turn to the case of $\Omega(n,2)$ and use Corollary \ref{c8}. Assume
first that $n=4$ and consider $N\in\Omega(4,2)$ corresponding to the
partition $(2,2)$. The only possibility for $G$ is then to be $B^2\times B^2$. 

Now let $n$ be arbitrary, and we
consider the case of $N\in\Omega(n,2,1)$ corresponding to the
partition $(n-1,1)$. Thus, in (\ref{b}) we have $p=2$, $n_1=n-1$,
$n_2=1$ and either $s=t=1$, or
$s=0,t=1$. Next, there are the following possibilities for a
hyperbolic Reinhardt domain $\widetilde G\subset\CC$:
\begin{eqnarray*}
&{}&\hbox{{\bf (a)}}\,\,\widetilde G=\{|z_n|<R\},\qquad 0<R<\infty;\\
&{}&\hbox{{\bf (b)}}\,\,\widetilde G=\{r<|z_n|<R\},\qquad 0<r<R\le\infty;\\
&{}&\hbox{{\bf (c)}}\,\,\widetilde G=\{0<|z_n|<R\},\qquad 0<R<\infty.
\end{eqnarray*}
Substituting {\bf (a)}, {\bf (b)}, {\bf (c)} into (\ref{b}), and excluding non-hyperbolic
domains, we produce {\bf (iv)}--{\bf (vi)} (cf. \cite{IK3}).

The theorem is proved. \qed
\smallskip\\

In Corollaries \ref{c2}, \ref{c5} and Theorem \ref{t9} we have described all
hyperbolic Reinhardt domains whose group of holomorphic automorphisms
has dimension $n^2$ or higher. We now turn to the case of dimensions
not exceeding $n^2-2$. To this end, we introduce the sets
\begin{itemize}
\item $\displaystyle \widetilde\Omega(n):=\{N\in\Omega(n):\, N\le n^2-2\}$,
\item $\displaystyle \widetilde\Omega(n,\ell):=
     \{N\in\Omega(n,\ell):\, N\le n^2-2\}$,
\item $\displaystyle \widetilde\Omega(n,\ell,q):=\{N\in\Omega(n,\ell,q):\, 
  N\le n^2-2\},\,\,\ell\ge 2$. 
\end{itemize}
The numbers
belonging to $\cup_{\ell=2}^n\Omega(n,\ell,0)$ 
we call ${\it compact}$ and the
numbers belonging to 
$\widetilde\Omega(n)\setminus\cup_{\ell=2}^n\Omega(n,\ell,0)$
we call {\it non-compact} (note that $\displaystyle \widetilde
\Omega(n,\ell,0)=\displaystyle \Omega(n,\ell,0)$ for $\ell\ge 2$).  

It is clear from (\ref{c}) that compact numbers
arise as the dimensions of the automorphism groups of domains $G$ for
which $\hbox{Aut}_0(G)=U(n_1)\times\dots\times U(n_{\ell})$, for some
partition $(n_1,\dots,n_{\ell})$ of $n$ with $\ell\ge 2$. For any such
partition one can construct many pairwise
non-equivalent hyperbolic (and even smoothly bounded) Reinhardt
domains for which the identity component of the automorphism groups
is $U(n_1)\times\dots\times U(n_{\ell})$. The construction is as follows. Choose
a set $Q\subset\RR^{\ell}_{+}$ (here $\RR^{\ell}_{+}:=\{(x_1,\dots,x_{\ell}):x_j\ge
0\}$) in such a way that
$$
D_Q:=\left\{(z^1,\dots,z^{\ell})\in\CC^n:(|z^1|,\dots,|z^{\ell}|)\in Q\right\}
$$
is a smoothly bounded domain in $\CC^n$ containing the
origin; here
$$
z^i:=\left(z_{n_1+\dots+n_{i-1}+1},\dots,z_{n_1+\dots+n_i}\right),
\qquad i=1,\dots,\ell.
$$
By
\cite{Su}, two bounded Reinhardt domains containing the origin are
holomorphically equivalent if and only if one is obtained from another
by dilations and permutations of coordinates. In \cite{FIK1}
we listed all smoothly bounded Reinhardt domains with non-compact automorphism
group. Since they all contain the origin, it is not difficult to
choose $D_Q$ such that it is not holomorphically equivalent to any of
the domains from \cite{FIK1} and thus ensure that $\hbox{Aut}_0(D_Q)$ is
compact and therefore, by (\ref{c}), is isomorphic to a product of unitary
groups. Further, $U(n_1)\times\dots\times
U(n_{\ell})\subset\hbox{Aut}_0(D_Q)$, and if necessary, one can vary $Q$
slightly to get $U(n_1)\times\dots\times U(n_{\ell})=\hbox{Aut}_0(D_Q)$. The
freedom in choosing $Q$ satisfying the above requirements is very
substantial, and one can produce (uncountably) many non-equivalent
domains that clearly cannot be classified in any reasonable way.    
Thus compact dimensions are
virtually ``unclassifiable'' and so are ``bad'' dimensions. On the other hand, it is clear from
(\ref{c}) that if $\hbox{dim}\,\hbox{Aut}(D)$ is non-compact, then
$\hbox{Aut}(D)$ is non-compact. In \cite{IK1} we
classified all bounded Reinhardt domains with $C^1$-smooth boundary
and non-compact automorphism group. Thus, at least in the
$C^1$-smooth and bounded
situation, non-compact dimensions are ``classifiable'' and so are termed
``good'' dimensions. 

Thus from now on we will restrict our considerations to bounded
Reinhardt domains with $C^1$-smooth boundary. It follows from
\cite{FIK1} and the discussion above that only numbers of the form (\ref{e}) with $m=0,1$ 
can be realizable as the dimensions of the automorphism groups of such
domains. On the other hand, for any such number it is not difficult to construct a
smoothly bounded domain whose automorphism group has dimension equal
to this number. Indeed, assume that in (\ref{e}) $m=1$ (the case $m=0$ has been
discussed above) and consider a number $N$ of the form
\begin{equation}
\sum_{i=1}^k n_i^2+2n_{i_0},\label{f}
\end{equation}
for some index $i_0$. Consider the domain
$$
D_{n_1,\dots,n_k}^{i_0}(s):=\left\{(z^1,\dots,z^{\ell})\in\CC^n: |z^{i_0}|^2+\sum_{i\ne i_0}|z^i|^{2s_i}<1\right\},
$$
where $s:=(s_1,\dots,s_{i_0-1},s_{i_0+1},\dots,s_k)$, $s_i\in\NN$, $s_i>1$, $s_i\ne s_j$ $(i\ne j)$. By an argument
similar to that in the proof of Theorem 1 in \cite{FIK2} one can now
explicitly determine $\hbox{Aut}_0(D_{n_1,\dots,n_k}^{i_0}(s))$. It then follows from the
explicit formulas that
$\hbox{dim}\,\hbox{Aut}(D_{n_1,\dots,n_k}^{i_0}(s))=N$. 

Therefore, in the situation of $C^1$-smooth bounded domains the sets of interest are
$$
\widetilde\Omega^{{\rm sb}}(n,\ell):=\Omega(n,\ell,0)\cup\widetilde\Omega(n,\ell,1),\quad\ell\ge2,\quad\hbox{and}\quad\widetilde\Omega^{{\rm sb}}(n):=\cup_{\ell=2}^n\widetilde\Omega^{{\rm
sb}}(n,\ell).
$$ 
The set $\widetilde\Omega^{{\rm sb}}(n)$ appears to
have an extremely 
irregular structure.

Let $C(n)$ and $H(n)$ denote, respectively, the sets of all compact and
non-compact dimensions from $\widetilde\Omega^{{\rm sb}}(n)$. Any
number from 
$C(n)$ does not exceed $n^2-2$ and  has the form 
$$
\sum_{i=1}^k n_i^2, 
$$
where $(n_1,\dots,n_k)$
is a partition of $n$. Any number from $H(n)$ does not exceed $n^2-2$
and has the form (\ref{f}). A number of the form (\ref{f}) can be written
as
$$
\sum_{1\le i\le k
      \atop 
         i\ne i_0} n_i^2+(n_{i_0}+1)^2-1.
$$
Therefore,
\begin{equation}
H(n)=(C(n+1)-1)\setminus(C(n)\cup\{n^2\}).\label{k}
\end{equation}
Let $c(n)$,
$h(n)$ denote the cardinalities of $C(n)$, $H(n)$
respectively. Clearly we have
\begin{eqnarray}
c(n)&\le&\frac{n(n-1)}{2}<\frac{n^2}{2},\label{kom}\\
h(n)&\le&\frac{n(n-1)}{2}<\frac{n^2}{2}.\label{stu}
\end{eqnarray}
Since for any partition $(n_1,\dots,n_k)$ of $n$, the $k+1$-tuple $(1,n_1,\dots,n_k)$
is a partition of $n+1$, we have
$$
C(n)\cup\{n^2\}\subset C(n+1)-1,
$$
and formula (\ref{k}) implies:
\begin{equation}
h(n)=c(n+1)-c(n)-1. \label{g}
\end{equation}
Therefore, determining $c(n)$ at the same time yields $h(n)$.

Finding $c(n)$, however, has proved to be a difficult task. In Section 2
below we list some of the results of numerical computations for up to $n=4065$.  In general, we have
$$
C(n)=\cup_{\ell=2}^n\Omega(n,\ell,0).
$$
One can obtain characterizations of the sets $\Omega(n,\ell,0)$,
$2<\ell\le n-1$, similar to that of $\Omega(n,2,0)$ from Proposition \ref{p7}. However, it is still not clear how one can calculate $c(n)$ by
using such characterizations since the sets $\Omega(n,\ell,0)$,
$\ell=2,\dots,n-1$, intersect each other in a chaotic manner.
Nevertheless, we have been able to determine the principal term in the
asymptotic behavior of $c(n)$.

\begin{theorem}\label{num} \sl We have
\begin{equation}
\lim_{n\rightarrow\infty}\frac{c(n)}{n^2}=\frac{1}{2}.\label{ner}
\end{equation}

\end{theorem}

\noindent{\bf Proof:} Let, for any integer $n\ge 0$,
\begin{equation}
\widehat C(n):=\{\sum_{i=1}^k n_i^2:\sum_{i=1}^k n_i=n\};\label{sets}
\end{equation}
and let $\hat c(n)$ denote the cardinality of the set $\widehat
C(n)$. Clearly, $\hat c(n)=c(n)+1$ for $n\ge 2$. 

First, we note that
\begin{equation}
\widehat C(n)=\bigcup_{0\le i<n}\biggl(\widehat C(i)+(n-i)^2\biggr).\label{uni}
\end{equation} 

Next, we define sequences of integers $f(0), f(1),\dots$, $g(0), g(1),
\dots$ and $k(1), k(2), \dots$ inductively as follows:
\begin{enumerate}
\item $f(0)=0$;
\item $k(1)=0$;
\item $f(n) = (n-k(n))^2 + f(k(n))$;
\item $g(n) = {1\over 2} (f(n)+n+4)$;
\item $k(n)=\sup \{ 0\leq \kappa < n : g(\kappa) \leq n \}$.
\end{enumerate}
We note that $k(n)$ is a non-decreasing sequence 
and $k(n)\ra\infty$ as $n\ra\infty$.

We shall now show by induction on $n$ that $f(n)$ has the same parity as
$n$. This clearly holds for $n=0$. For arbitrary $n$ we have
$$
f(n)-n\equiv (n^2-n)+k^2(n)+f(k(n))\,\hbox{(mod $2$)}\equiv
k^2(n)+f(k(n))\,\hbox{(mod $2$)}.
$$
But, by induction, $f(k(n))\equiv k(n)\,\hbox{(mod
$2$)}$ and thus $f(n)-n\equiv 0\,\hbox{(mod $2$)}$. 

We shall now show by induction on $n$ that
$\{n,n+2,\dots,f(n)\}\subset\widehat C(n)$ (we will see below that each of
the sets  $\{n,n+2,\dots,f(n)\}$ is non-empty). Certainly the claim
is true for $n=0$. For the general case, let $j\in \{n,n+2,\dots,
f(n)\}$.
\smallskip\\

\noindent {\bf (a)} If $j\geq (n-k(n))^2 + k(n)$, then clearly $j\in \{k(n),k(n)+2,\dots,
f(k(n))\}+(n-k(n))^2$. By induction,
$\{k(n),k(n)+2,\dots,f(k(n))\}\subset\widehat C(k(n))$. Then it
follows from (\ref{uni}) that $j\in\widehat C(n)$. 

\noindent {\bf (b)} Suppose now that
$j<(n-k(n))^2+k(n)$. We define $k:=\sup\{\kappa:k(n)\leq
\kappa<n\,\,\mbox{and}\,\, j<
(n-\kappa)^2+\kappa\}$. Since $j\geq n$, we have $k<n-1$. It now follows from
the definitions of $k$ and $f(n)$ and the fact that $g(k+1)>n$ that
$j\in\{k+1,k+3,\dots,f(k+1)\}+(n-(k+1))^2$; as above, by induction
and (\ref{uni}) we obtain that $j\in C(n)$.

Next, we claim that $f(n)\geq 2n$ for $n\geq 4$. This is verified by
explicit calculation for small $n$; in particular it is true for
$4\leq n\leq 18$. For general $n$, we have $f(n)\geq (n-k(n))^2$.
Since $n\geq g(k(n))={1\over 2}(f(k(n)+k(n)+4)$ and $k(n)>4$ for $n\ge 18$
(one can check that $k(18)=7$ and use the fact that $k(n)$ is non-decreasing),
we have by induction that $n\geq {3\over 2}k(n)$ and so
\begin{equation}
f(n)\geq (n-k(n))^2\geq {1\over 9}n^2\geq 2n, \qquad \hbox{for $n\geq
18$}.\label{dop}
\end{equation}

The inequalities (\ref{dop}) also imply that $f(n)/n\ra \infty$ as $n\ra\infty$. Since
$k(n)\ra \infty$ as $n\ra\infty$ and $n\geq
g(k(n))\geq {1\over 2}f(k(n))$, we find that $n/k(n)\ra \infty$ as
$n\ra\infty$, hence
\begin{equation}
{f(n)\over n^2}\geq\left({ n-k(n)\over n}\right)^2\ra1,\qquad \hbox{as}\,\, n\ra\infty.\label{oc}
\end{equation}
Since the set $\{n,n+2,\ldots,f(n)\}$ is contained in $\widehat C(n)$ for
all $n$, we have
$$
c(n)=\hat c(n)-1\ge \frac{f(n)-n}{2}.
$$
Combining this inequality with (\ref{kom}) and (\ref{oc}) we obtain (\ref{ner}).

The theorem is proved.\qed
\smallskip\\

We shall now find a lower bound for $\liminf_{n\ra\infty}
h(n)/n$. First, we need the following technical lemma.

\begin{lemma} \label{largest} \sl Let $\widehat C(n)$ be the sets
defined in (\ref{sets}), and suppose that $n\ge 7$. If $N\in\widehat C(n)$ and $N>3n^2/4$, then
$N\in\widehat C(i)+(n-i)^2$ for some $i<n/2$.
\end{lemma}

\noindent{\bf Proof:} First, we shall show that if
$\lambda=(n_1,\dots,n_k)$ is a partition of $n$ such that $n_j\le n/2$
for all $j$, then
\begin{equation}
\sum_{j=1}^k n_j^2\le \frac{3n^2}{4}.\label{3/4}
\end{equation}
Let $X(n)$ denote the maximal value of the sum $\sum_{j=1}^k n_j^2$ over
all such partitions $\lambda$. Choose a partition
$\lambda'=(n_1',\dots,n_m')$ such that $X(n)=\sum_{j=1}^m
{n_j'}^2$. We claim that $m\le 3$. Indeed, suppose that $m\ge
4$. Consider the set $S$ of all indices $j\in\{1,\dots,m\}$ such that 
$$
n_j'+1\le\frac{n}{2}.
$$
Clearly, if $n\ge 7$, the set $S$ contains at least two elements. Let
$j_1,j_2\in S$, $j_1\ne j_2$, and suppose that $n_{j_1}\ge n_{j_2}$. Consider a partition
$\tilde\lambda=(\tilde n_1,\dots,\tilde n_m)$ defined as follows
\begin{eqnarray*}
\tilde n_j&:=&n_j,\qquad\hbox{if $j\ne j_1,j_2$},\\
\tilde n_{j_1}&:=&n_{j_1}+1,\\
\tilde n_{j_2}&:=&n_{j_2}-1.
\end{eqnarray*}
Then
$$
\sum_{j=1}^m \tilde n_j^2=X(n)+2(n_{j_1}-n_{j_2})+2>X(n),
$$
which contradicts the definition of $X(n)$. Therefore $m\le 3$ which implies that
$X(n)\le 3n^2/4$, and thus (\ref{3/4}) holds.

Thus, the number $N$ can only be realizable by a partition for which
at least one entry is bigger than $n/2$, and the lemma is proved.\qed
\smallskip\\

We are now ready to prove the following theorem. 

\begin{theorem}\label{numh} \sl We have
$$
\liminf_{n\rightarrow\infty}\frac{h(n)}{n}\ge 1.
$$

\end{theorem}

\noindent {\bf Proof:} Let $\widehat C(n)$ be the sets defined in
(\ref{sets}). Define
$$
\widehat H(n):=(\widehat C(n+1)-1)\setminus \widehat C(n),
$$
for $n=0,1,\dots$ Let $\hat h(n)$ denote the cardinality of $\widehat
H(n)$. Since $\widehat C(n)\subset \widehat C(n+1)-1$, we have
\begin{equation}
\hat h(n)=\hat c(n+1)-\hat c(n),\label{incl}
\end{equation}
for all $n$, where $\hat c(n)$ is the cardinality of $\widehat C(n)$.

We now fix $n$ and choose a number $k\in\NN$ such that
\begin{equation}
\frac{k^2+3k+1}{2}\le n<\frac{(k+1)^2+3(k+1)+1}{2}.\label{ok}
\end{equation}
It follows from (\ref{uni}) that
\begin{eqnarray}
\widehat H(n)&=&\cup_{0\le i<n+1}\left(\widehat C(i)+(n+1-i)^2-1\right)\setminus \nonumber\\
&{}&\cup_{0\le
i<n+1}\left(\widehat C(i-1)+(n+1-i)^2\right),\label{cru}
\end{eqnarray}
where we set $\widehat C(-1):=\emptyset$. 

We observe that the sets
$\widehat C(i)+(n+1-i)^2-1$ for $i\le k$ lie strictly above $3(n+1)^2/4$, if $n$ is
sufficiently large. To prove this, we show that
$(n+1-i)^2-1>3(n+1)^2/4$ for $i\le k$. Indeed, it
follows from the first inequality in (\ref{ok}) that for large $n$ one has
$i\le k\le 3\sqrt{n}/2$. Therefore
$$
(n+1-i)^2-1\ge (n+1)^2-3\sqrt{n}(n+1)+\frac{9n}{4}-1.
$$
The expression on the right-hand side is clearly bigger than
$3(n+1)^2/4$ for large $n$.
\smallskip\\

\noindent{\bf Claim.} The sets $\widehat C(i)+(n+1-i)^2-1$ and $\widehat
C(j)+(n+1-j)^2-1$ do not intersect if $0\le i\le k$, $0\le j<(n+1)/2$,
$i\ne j$.
\smallskip\\

To prove the claim we note that $\widehat C(\ell)+(n+1-\ell)^2-1\subset
[\ell,\ell^2]+(n+1-\ell)^2-1$ for any $\ell$. Let $0\le i,i+1\le
k$. It then follows from the first inequality in (\ref{ok}) that
$(i+1)^2+(n+1-(i+1))^2-1<i+(n+1-i)^2-1$. This inequality also holds true if
$i=k$. Since $\ell^2+(n+1-\ell)^2-1$ is a decreasing function of
$\ell$ for $\ell<(n+1)/2$, the claim follows.

Further, Lemma \ref{largest} and formula (\ref{cru}) imply
$$
\widehat H(n)\supset\cup_{0\le i\le k}\left(\widehat
C(i)+(n+1-i)^2-1)\setminus \widehat C(i-1)+(n+1-i)^2\right),
$$
and therefore by (\ref{incl}) 
\begin{equation}
\hat h(n)\ge\sum_{i=1}^k \hat h(i)=\hat c(k).\label{h}
\end{equation}

Let $0\le\beta<1$ and $n\ge (k_0^2+3k_0+1)/2$ where $k_0$ is chosen to
satisfy the following
\begin{eqnarray*}
&{}&\frac{k_0^2}{k_0^2+5k_0+5}\ge\sqrt{\beta},\\
&{}&\hat c(k)\ge\frac{\sqrt{\beta}k^2}{2},\qquad\hbox{for all $k\ge
k_0$,}
\end{eqnarray*}
(the second inequality can be satisfied by Theorem \ref{num}).
From (\ref{h}) and the second inequality in (\ref{ok}) we now have
\begin{eqnarray*}
&{}&\frac{\hat h(n)}{n}\ge
\frac{\sqrt{\beta}k^2}{(k+1)^2+3(k+1)+1}=\\
&{}&\frac{\sqrt{\beta}k^2}{k^2+5k+5}\ge
\frac{\sqrt{\beta}k_0^2}{k_0^2+5k_0+5}\ge \beta.
\end{eqnarray*}
Since  $h(n)=\hat h(n)-1$, the theorem follows.\qed
\smallskip\\

\noindent{\bf Remark.} It is plausible---especially in view of our
numerical computations (see Section 2)--- that in fact one has
$\lim_{n\rightarrow\infty}h(n)/n=1$. It is straightforward to show that
$\limsup_{n\rightarrow\infty}h(n)/(nk(n))\le 1$, where $k(n)$ are the
numbers defined in the proof of Theorem \ref{num}. Note that $k(n)$
grows much more slowly than $n$ (cf. (\ref{stu})).
 
We will now look at the numbers $c(n)$ from the point of view of the
theory of partitions (see e.g. \cite{A} for a general exposition of
this theory). Let $\lambda=(n_1,\dots,n_k)$ be a partition of $n$ and
assume that $n_1\ge n_2\ge\dots\ge n_k$. The partition $\lambda$ can
be pictured by utilizing its {\it Young diagram} $Y(\lambda)$:

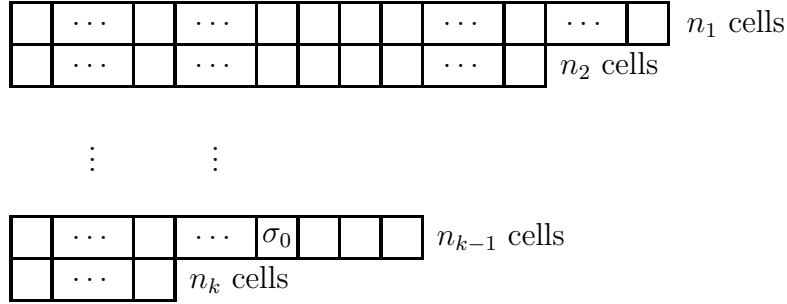
\begin{picture}(300,200)(-30,-30)
\thicklines
\put(0,0){\framebox(15,15)}
\put(15.5,0){\framebox(30,15){\dots}}
\put(46.3,0){\framebox(15,15)}
\put(66,0){\makebox(40,15){$n_k$ cells}}

\put(0,15.8){\framebox(15,15)}
\put(15.5,15.8){\framebox(30,15){\dots}}
\put(46.3,15.8){\framebox(15,15)}
\put(62,15.8){\framebox(30,15){\dots}}
\put(92.8,15.8){\framebox(15,15){$\sigma_0$}}
\put(108.5,15.8){\framebox(15,15)}
\put(124.2,15.8){\framebox(15,15)}
\put(140,15.8){\framebox(15,15)}
\put(165,15.8){\makebox(40,15){$n_{k-1}$ cells}}

\put(15.5,21){\makebox(30,70){\vdots}}
\put(62,21){\makebox(30,70){\vdots}}

\put(0,81){\framebox(15,15)}
\put(15.5,81){\framebox(30,15){\dots}}
\put(46.3,81){\framebox(15,15)}
\put(62,81){\framebox(30,15){\dots}}
\put(92.8,81){\framebox(15,15)}
\put(108.5,81){\framebox(15,15)}
\put(124.2,81){\framebox(15,15)}
\put(140,81){\framebox(15,15)}
\put(155.8,81){\framebox(30,15){\dots}}
\put(186.5,81){\framebox(15,15)}
\put(206,81){\makebox(40,15){$n_2$ cells}}

\put(0,96.8){\framebox(15,15)}
\put(15.5,96.8){\framebox(30,15){\dots}}
\put(46.3,96.8){\framebox(15,15)}
\put(62,96.8){\framebox(30,15){\dots}}
\put(92.8,96.8){\framebox(15,15)}
\put(108.5,96.8){\framebox(15,15)}
\put(124.2,96.8){\framebox(15,15)}
\put(140,96.8){\framebox(15,15)}
\put(155.8,96.8){\framebox(30,15){\dots}}
\put(186.5,96.8){\framebox(15,15)}
\put(202.3,96.8){\framebox(30,15){\dots}}
\put(233,96.8){\framebox(15,15)}
\put(254,96.8){\makebox(40,15){$n_1$ cells}}

\end{picture}

For any cell $\sigma\in Y(\lambda)$ one can define the {\it arm} of
$\sigma$ (denote it by $a(\sigma)$) as the number of cells to the right of $\sigma$. In the
diagram above, $a(\sigma_0)=3$. Clearly, the sum of all arms of the
cells in the $j^{\hbox{\small th}}$ row of $Y(\lambda)$ is equal to
$$
\sum_{m=0}^{n_j-1}m=\frac{n_j(n_j-1)}{2},
$$
and therefore the total number $A(\lambda)$ of arms in $Y(\lambda)$ is
$$
A(\lambda)=\sum_{j=1}^k
\frac{n_j(n_j-1)}{2}=\frac{1}{2}\left(\sum_{j=1}^k n_j^2-n\right).
$$
Therefore, $c(n)$ is equal to the number of distinct values that $A(\lambda)$
takes over all partitions $\lambda$ of $n$.

It would be very interesting to know whether the above
interpretation of $c(n)$ can be used to get more information on $c(n)$
and $h(n)$
by applying
techniques from the theory of partition.

\section{Numerical Computations for $c(n)$}
\setcounter{equation}{0}

Equation (\ref{uni}) can be used to form an efficient algorithm
for the calculation of $\widehat C(n)=C(n)\cup\{n^2\}$. Once the sets $\hat C(k)$ are known for $1\leq
k\leq n-1$, the set $\hat C(n)$ is formed by translating each $\hat C(k)$ and
taking the union. If $\hat C(k)$ is stored as an array of $O(k^2)$
elements, this process will take $O(\sum_{k<n} k^2)=O(n^3)$
time. Therefore to calculate all the sets $\hat C(k)$ for $1\leq k\leq n$
will take $\sum_{k\leq n}O(k^3)=O(n^4)$ time and $O(n^3)$ memory.

A C program implementing this algorithm has been written and used to
calculate $\widehat C(n)$ for $n\leq 2400$. As previously noted, the source code of this program is
available on the World Wide Web. The program stores
each element of $\hat C(n)$ as a single bit, which makes the program
more difficult to understand, but increases its speed and decreases its 
memory consumption considerably.
To give a sense of the constants
involved in the $O(n^4)$ time and $O(n^3)$ memory estimates, the
program takes 70 seconds to calculate to $n=1000$ on a Sun Ultra 10, and
uses around 20 megabytes of memory.

Below we list some results of numerical computations for
$n$ up to $4065$. They were obtained by combining the above-mentioned
C programming with network programming. We compute $c(n)=\hat c(n)-1$ and
then by applying (\ref{g}) calculate the numbers $h(n)$. 
Further, we compare the growth of $c(n)$ with $n^2$ and the growth of
$h(n)$
with $n$.
\medskip\\

\begin{center}
\begin{tabular}{|r|r|l|c|l|}
\hline
\multicolumn{1}{|c|}{$n$}&
\multicolumn{1}{c|}{$c(n)$}&
\multicolumn{1}{c|}{$c(n)/n^2$}&
\multicolumn{1}{c|}{$h(n)$}&
\multicolumn{1}{c|}{$h(n)/n$}
\\ \hline
20 & 117 & 0.2425 & 11 & 0.5500
\\ \hline
40 & 537 & 0.3356 & 31 & 0.7750
\\ \hline
60 & 1294 & 0.3595 & 47 & 0.7833
\\ \hline
80 & 2403 & 0.3754 & 62 & 0.7750
\\ \hline
100 & 3880 & 0.3880 & 81 & 0.8100
\\ \hline
200 & 16785 & 0.4196 & 176 & 0.8800
\\ \hline
400 & 70922 & 0.4432 & 365 & 0.9125
\\ \hline
600 & 163415 & 0.4539 & 559 & 0.9316
\\ \hline
800 & 294630 & 0.4603 & 753 & 0.9412
\\ \hline
1000 & 464692 & 0.4646 & 949 & 0.9490
\\ \hline
1500 & 1060777 & 0.4714 & 1439 & 0.9593
\\ \hline
2000 & 1901854 & 0.4754  & 1925  & 0.9625
\\ \hline
2500 & 2988578 & 0.4781  & 2423 & 0.9691
\\ \hline
3000 & 4321549  & 0.4801 & 2908  & 0.9693
\\ \hline
3500 & 5901013 & 0.4817 & 3410 & 0.9742
\\ \hline
4000 & 7731988 & 0.4832 & 3466 & 0.8665
\\ \hline
4065 & 7982961 & 0.4831 & 3962 & 0.9746
\\ \hline
\end{tabular}
\medskip\\
\end{center}

\section{Concluding Remarks}
\setcounter{equation}{0}

In this paper we have endeavored to correlate the dimension
of the automorphism group of a domain in $\CC^n$ with
the geometric characteristics of the domain.  By restricting
attention to Reinhardt domains, we have been able to exploit
the structure theorem of Kruzhilin and to come up
with (at least in some cases) rather specific conclusions.

We hope that the information obtained here in the Reinhardt case
points the way toward what ought to be true for more general classes
of domains.  In particular, we have identified certain dimensions for
the automorphism group that we call ``good'' and certain others that
we call ``bad''.  The former are dimensions in which the automorphism
groups are always non-compact; certainly the extant literature (see \cite{IK2})
suggests that the Levi geometry of a boundary accumulation point
gives us a chance of classifying such domains. The latter 
are dimensions for which there exist compact automorphism groups.
We are able to determine explicitly that in these cases a holomorphic
classification does not exist. 

We have been able to find the principal term in the asymptotics
of the number of ``bad''
automorphism group dimensions as well as to bound from
below the principle term in the asymptotics of the number of ``good'' ones.  In particular, the ``good''
dimensions are fairly robust as $n \ra \infty$; this is positive information.  
We also utilize a computer
to count these numbers in complex dimension
$n$ for values of $n$ up to $4065$. 

It is clear that further effort is needed to show that
$h(n)/n\ra\infty$ as well as to determine the forms of
the error terms $c(n)-n^2/2$ and $h(n)-n$ that our numerical
computations can only suggest. We plan
to develop these ideas in future papers. In particular, we wish to
extend the work to domains more general than Reinhardt domains.

\bigskip

{\obeylines
School of Mathematical Sciences
The Australian National University 
Canberra, ACT 0200
AUSTRALIA 
E-mail address: James.Gifford@anu.edu.au
\smallskip

Centre for Mathematics and Its Applications 
The Australian National University 
Canberra, ACT 0200
AUSTRALIA 
E-mail address: Alexander.Isaev@anu.edu.au
\smallskip

Department of Mathematics
Washington University, St.Louis, MO 63130
USA 
E-mail address: sk@math.wustl.edu}

\end{document}